\documentclass{llncs}
\usepackage[]{amsmath}
\usepackage[]{amssymb}
\usepackage{graphicx}

\newcommand{\abs}[1]{{\left|#1\right|}}
\newcommand{\set}[1]{{\left\{#1\right\}}}
\newcommand{\st}{\colon}

\newcommand{\RP}{\mathop\mathsf{RP}}
\newcommand{\PR}{\mathop\mathsf{PR}}
\newcommand{\RPPR}{\mathop\mathsf{RPPR}}
\newcommand{\ANY}{\mathop\mathsf{ANY}}

\newcommand{\ord}{{\mathop{\mathrm{ord}}\nolimits}}

\newcommand{\Nfp}{F}

\newcommand{\sigzero}{d}
\newcommand{\sigone}{\sigma}

\spnewtheorem{heuristic}[theorem]{Heuristic}{\bfseries}{\itshape}

\begin{document}

\title{Distribution of the Error in Estimated Numbers of Fixed Points
 of the Discrete Logarithm} 
\author{Joshua Holden}
\institute{Department of Mathematics, Rose-Hulman Institute of
Technology, Terre Haute, IN, 47803-3999, USA,
\email{holden@rose-hulman.edu}}

\maketitle

\begin{abstract}
Brizolis asked the question: does every prime $p$ have a pair $(g,h)$
such that $h$ is a fixed point for the discrete logarithm with base
$g$?  This author and Pieter Moree, building on work of Zhang, Cobeli,
and Zaharescu, gave heuristics for estimating the number of such pairs
and proved bounds on the error in the estimates.  These bounds are not
descriptive of the true situation, however, and this paper is a first
attempt to collect and analyze some data on the distribution of the
actual error in the estimates.
\end{abstract}

\section{Introduction}

Paragraph~F9 of~\cite{UPINT} includes the following problem,
attributed to Brizolis: given a prime $p>3$, is there always a
pair $(g,h)$ such that $g$ is a primitive root of $p$, $1 \leq h
\leq p-1$, and
\begin{equation} \label{fp}
g^{h} \equiv h  \mod{p} \enspace ?
\end{equation}
In other words, is there always a primitive root $g$ such that the
discrete logarithm $\log_{g}$ has a fixed point?  This question has now
been settled affirmatively by Campbell and Pomerance
in~\cite{CampbellThesis}.  The answer relies on an estimate for the
number $N(p)$ of pairs $(g,h)$ which satisfy the equation, have $g$ is
primitive root, and also have $h$ a primitive root which thus must be
relatively prime to $p-1$.  This result seems to have been discovered
and proved by Zhang in~\cite{Zhang} and later, independently, by Cobeli
and Zaharescu in~\cite{CZ}.

In~\cite{HM2004} and~\cite{HM2005}, Pieter Moree and this author 
applied the same methods to estimate the number of solutions 
to~\eqref{fp} given no conditions on $g$ and $h$.  Unfortunately, the 
error term involved in this estimate was completely unsatisfactory.  
It was also shown in~\cite{HM2005} that for a positive proportion of 
primes a better error estimate can be obtained, and it was conjectured 
that one could do even better.  The object of this note is to collect 
and analyze some data on the distribution of the actual error in these 
estimates.

The idea of repeatedly applying the function $x \mapsto g^{x}
\bmod{p}$ is used in the famous cryptographically secure pseudorandom
bit generator of Blum and Micali.  (\cite{Blum-Micali}; see
also~\cite{Patel-Sundaram} and~\cite{Gennaro}, among others, for
further developments.)  If one could predict that a pseudorandom
generator was going to fall into a fixed point or cycle of small
length, this would obviously be detrimental to cryptographic security.
We hope that the investigation of the cycle structure of the discrete 
logarithm will therefore eventually be of some use to those interested 
in the field of cryptography.

Using the same notation as in the previously cited papers, we will refer to an
integer which is a primitive root modulo $p$ as $\PR$ and an integer
which is relatively prime to $p-1$ as $\RP$.  An integer which is both
will be referred to as $\RPPR$ and one which has no restrictions will
be referred to as $\ANY$.  
% In some instances, $\bullet$ will be used to stand for any one of
% these four conditions.

All integers will be taken to be between $1$ and $p-1$, inclusive,
unless stated otherwise.  If $N(p)$ is, as above, the number of
solutions to~\eqref{fp} such that $g$ is a primitive root and $h$ is a
primitive root which is relatively prime to $p-1$, then we will say
$N(p)=\Nfp_{g \PR, h \RPPR}(p),$ and similarly for other
conditions.  
% If $\ord_{p}(g)=\ord_{p}(h)$, we say that $g \ORD h$.
We will use $\sigzero(n)$ for the number of divisors of
$n$ and $\sigone(n)$ for the sum of the divisors of $n$.  All other
notations should be fairly standard.

\section{Heuristics, Conjectures, and Previous Results}

The fundamental observation at the heart of the estimation of
$\Nfp_{g \PR, h \RPPR}(p)$ is that if $h$ is a primitive root
modulo $p$ which is also relatively prime to $p-1$, then there is a
unique primitive root $g$ satisfying~\eqref{fp}, namely $g =
h^{\overline{h}}$ reduced modulo $p$, where $\overline{h}$ denotes the
inverse of $h$ modulo $p-1$ throughout this note.  Thus to estimate
$N(p)$, we only need to count the number of such $h$; $g$ no longer
has to be considered.  We observe that there are $\phi(p-1)$
possibilities for $h$ which are relatively prime to $p-1$, and we
would expect each of them to be a primitive root with probability
$\phi(p-1)/(p-1)$.  This heuristic uses the assumption that the
condition of being a primitive root is in some sense ``independent''
of the condition of being relatively prime.

We will actually need the following slightly more general heuristic:

\begin{heuristic}[Heuristic~2.6 of~\cite{HM2004}] \label{indepheur}
    The order of $x$ modulo $p$ is independent of the
greatest common divisor of $x$ and $p-1$,  in the sense that for all $p$, 
and all divisors $e$ and $f$ of $p-1$,
    \begin{multline*}
    \frac{1}{p-1}\#\set{x \in \set{1, \ldots, p-1} \st \gcd(x,p-1)=e, \quad
\ord_{p}(x)=\frac{p-1}{f} } \\
    \begin{aligned}
    &\approx \frac{1}{p-1}  \#\set{x \in \set{1, \ldots, p-1} \st \gcd(x,p-1)=e} \\
    &   \quad \times
    \frac{1}{p-1} \#\set{x \in \set{1, \ldots, p-1} \st \ord_{p}(x)=\frac{p-1}{f}}.
    \end{aligned}
    \end{multline*}
\end{heuristic}

The following lemma makes this heuristic rigorous; it was stated and proved 
in~\cite{HM2004} using the ideas in~\cite{CZ}.
\begin{lemma}[Lemma~2.7 of~\cite{HM2004}] \label{ordgcdlemma}
    Let $e$ and $f$ be divisors of
$p-1$, and $N$ a multiple of $p-1$.  Let
$\mathcal{P}=\set{1, \ldots, N}$ and
\[
\mathcal{P}' = \left\{ x \in \mathcal{P} \st \gcd(x,p-1)=e, \quad
\ord_{p}(x)=\frac{p-1}{f} \right\}.
\]
Then
\begin{multline*}
\begin{aligned}
    \abs{\#\mathcal{P}' - \frac{N}{(p-1)^2} \phi\left(\frac{p-1}{f}\right)
\phi\left(\frac{p-1}{e}\right)}
&\leq
 \sigzero \left(\frac{p-1}{f}\right)
\sigzero\left(\frac{p-1}{e}\right) \sqrt{p}(1+\ln p) \\
&\leq
 \sigzero\left(p-1\right)^{2}
 \sqrt{p}(1+\ln p).
\end{aligned}
\end{multline*}
\end{lemma}

Using this lemma with $e=f=1$ it is straightforward to prove Cobeli and
Zaharescu's version of Zhang's result.

\begin{theorem}[Theorem~1 of~\cite{CZ}] \label{cz1thm}
\[
\abs{\Nfp_{g \PR, h \RPPR}(p)  - \frac{\phi(p-1)^{2}}{p-1}} \leq
\sigzero(p-1)^{2}\sqrt{p}(1+\ln p).
\]
\end{theorem}

For the situation with no conditions on $g$ and $h$, we see
that~\eqref{fp} can be solved exactly when $\gcd(h, p-1)=e$ and $h$ is
a $e$-th power modulo $p$, and in fact there are exactly $e$ such
solutions.  Thus
\begin{equation} \label{NeTeqn}
\Nfp_{g \ANY, h \ANY}(p) = \sum_{e \mid p-1} e\ T(e,p).
\end{equation}
where
\[
T(e,p) =
\#\left\{ h \in
\mathcal{P}\left(1,1,{p-1}\right)^{(e)} \st \gcd(h,p-1)=e \right\}.
\]

Applying the lemma with $e=f=1$ gives
\begin{proposition}[Proposition~4.2 of~\cite{HM2004}] \label{anyanyprop}
Let $ e \mid p-1$.  Then
    \begin{enumerate}
	\item \label{Tpart} $\displaystyle \abs{T(e,p) - \frac{1}{e} \phi\left(\frac{p-1}{e} \right)}
	      \leq \sigzero \left(\frac{p-1}{e}\right) \sqrt{p}(1+\ln p).$
	\item \label{T1part} $\displaystyle T(1,p) = \phi(p-1).$
	\item \label{Tlargepart} $\displaystyle T(p-1,p) 
	= T\left(\frac{p-1}{2}, p\right) = 0.$
	\item \label{Ttrivpart} $\displaystyle 0 \leq T(e,p) \leq 
	\phi\left(\frac{p-1}{e}\right).$
	\item \label{anyanypart} 
	    \begin{multline*}
    \abs{\Nfp_{g \ANY, h \ANY}(p)  - (p-3)} \\
    \leq \sigzero(p-1)
    \left(\sigone(p-1)-\frac{3}{2}(p-1)\right) \sqrt{p}(1+\ln p).
    \end{multline*}

    \end{enumerate}
\end{proposition}

Unfortunately, the ``error'' term in Part~(\ref{anyanypart}) will be
larger than the main term for infinitely many $p$.  Using the deep 
result of Fouvry (see, e.g., 
\cite{Fouvry}) that $\gg x/\ln x$ primes $p \leq x$ are such that 
$p-1$ has a prime factor larger than $p^{0.6687}$, it was proved that:

\begin{theorem}[Theorem~4.8 of~\cite{HM2004}] \label{pieterA}
    There are $\gg x/\ln x$ primes $p \leq x$ such that
	\[
    \Nfp_{g \ANY, h \ANY}(p)  = (p-1) + O\left( 
    p^{5/6}\right).
    \]
    More specifically, there are $\gg x/\ln x$ primes $p \leq x$ such that
    \[
       \abs{\Nfp_{g \ANY, h \ANY}(p) - (p-1)} \leq
       p^{0.8313} \sigzero(p-1)^{2} (2 + \ln p).       
    \]
\end{theorem}

It was also noted in~\cite{HM2004} that if Fouvry's assertion holds
true with 0.6687 replaced by some larger $\theta$ (up to $\theta =
3/4$), then in Theorem~\ref{pieterA} the exponents $5/6$ and $0.8313$
can be replaced by $3/2-\theta+\delta$ and $3/2-\theta$ for any $\delta>0$.

On the other hand, we also expect that for many primes the error term 
cannot be set too small.  According to Heuristic~\ref{indepheur}, we can model
$T(e,p)$ using a set of independent random variables $X_{1}, \ldots,
X_{p-1}$ such that
\[
X_{h} = \begin{cases}
    \gcd(h,p-1) & \text{with probability $\frac{1}{\gcd(h,p-1)}$}; \\
    0 & \text{otherwise.}
\end{cases}
\]
Then the heuristic suggests that $\Nfp_{g \ANY, h \ANY}(p)$ 
is approximately equal to the expected value of $X_{1} + \cdots + 
X_{p-1}$, which is clearly $p-1$.  On the other hand, the variance 
$\sigma^{2}$ is the 
expected value of
\[
\left(\sum_{h=1}^{p-1} X_{h}-(p-1) \right)^{2}.
\]
Note that the expected value of $X_{h}X_{j}$ is $\gcd(h, p-1)$ if 
$h=j$ and $1$ otherwise.  Using this, an easy computation shows that
\begin{align*}
    \sigma^{2} &= \sum_{h=1}^{p-1} \gcd(h,p-1) - (p-1) 
     = \sum_{d \mid p-1} d \ \phi\left(\frac{p-1}{d}\right) - (p-1).
\end{align*}
In particular, the standard deviation $\sigma$ is less than
$p^{1/2+\epsilon}$ for every $\epsilon>0$ (for sufficiently large $p$).
Thus we have the following:

\begin{conjecture}[Conjecture~3.6 of~\cite{HM2004}] \label{varconj}
There are $o(x/\ln x)$ primes $p \leq x$ for which
\[
\abs{N_{\eqref{fp},  g \ANY, h \ANY}(p) - (p-1)} > p^{1/2+\epsilon}
\]
for every $\epsilon>0$.
\end{conjecture}

\section{Data and Analysis}
\begin{figure}[h!] 
	\caption{Plot of values of $\delta \geq 0$ for $3 \leq p \leq 15413$} 
	\label{posplot}
\begin{center}
\includegraphics[height=3in]{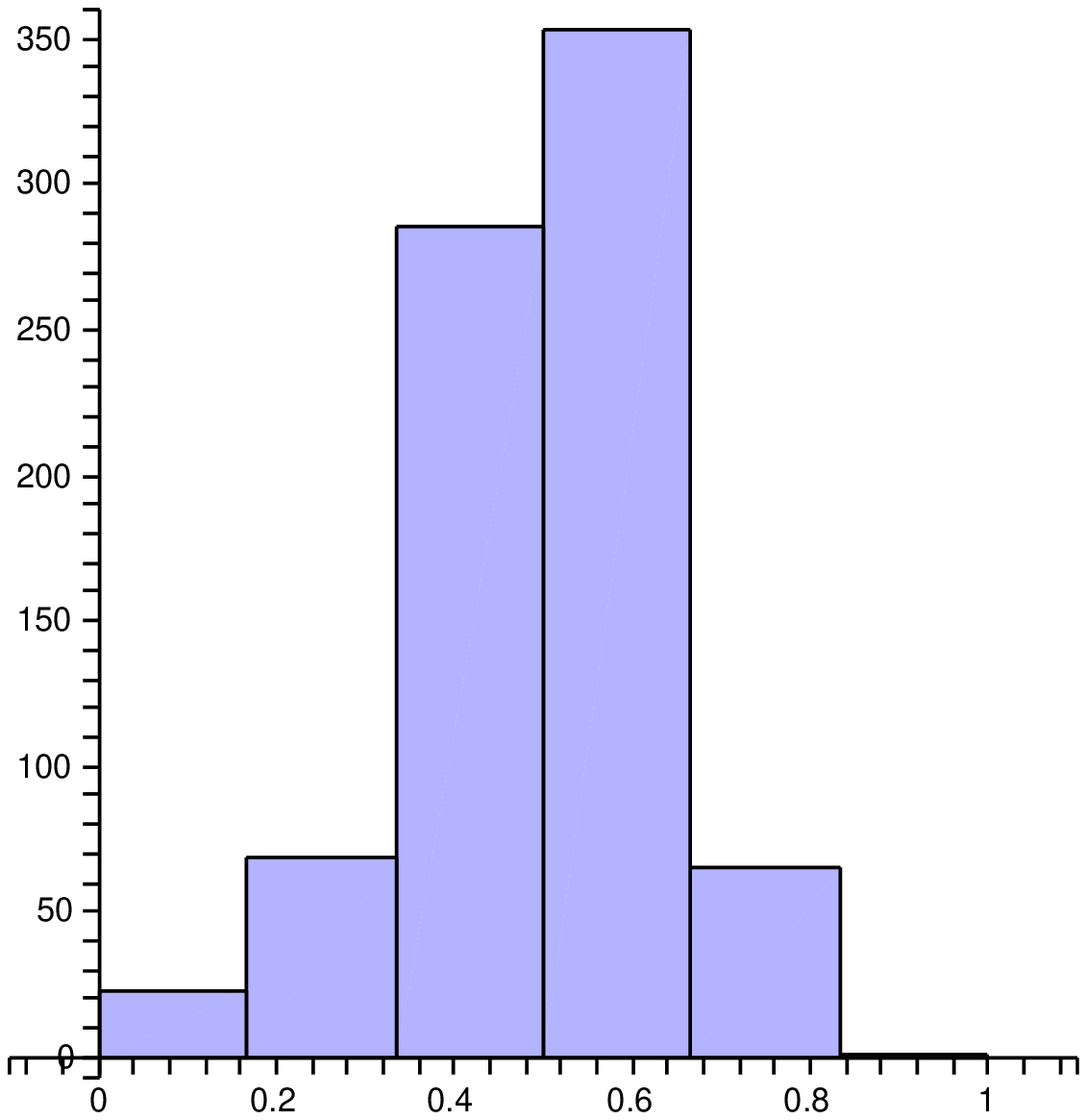}
\end{center}
\end{figure}

\begin{figure}[h!] 
	\caption{Plot of values of $\delta < 0$ for $3 \leq p \leq 15413$} 
	\label{negplot}
\begin{center}
\includegraphics[height=3in]{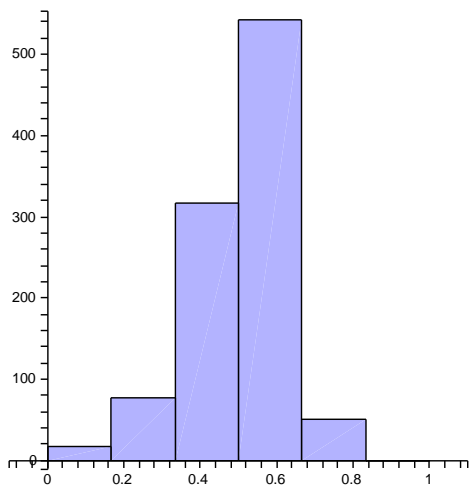}
\end{center}
\end{figure}

Since a factor of the form $p^{\alpha}$ dominates all of the proven
and conjectured bounds on the error given above, we decided to collect
data on the values of $\delta = {N_{\eqref{fp}, g \ANY, h \ANY}(p) -
(p-3)}$ for the first 1800 primes (3 through 15413).  The data was
then tallied based on the value of $\log_{p} \abs{\delta}$.
Table~\ref{posdelta} and Figure~\ref{posplot} give the data for
$\delta \geq 0$, while Table~\ref{negdelta} and Figure~\ref{negplot}
give the data for $\delta<0$.  The case $\delta=0$ did not actually
occur in this sample.  Likewise, there were no cases where
$\abs{delta} > p$, although this is certainly not ruled out for
$\delta > 0$.

\begin{table}[h!] 
	\caption{Values of $\delta \geq 0$ for $3 \leq p \leq 15413$} 
	\label{posdelta}
\begin{center}
	\begin{tabular}{||c|c|c|c|c|c|c||c||}
	    \hline
		$\log_{p} \abs{\delta}$ & 0--${1/6}$ & ${1/6}$--${1/3}$ & 
		${1/3}$--${1/2}$ & ${1/2}$--${2/3}$ & ${2/3}$--${5/6}$ &
		${5/6}$--$1$  & total	\\
		\hline
		$\#$ of $p$ &  23 & 69 & 285 & 353 & 65 & 1 & 796\\
		\hline
	\end{tabular}
\end{center}
\end{table}

\begin{table}[h!] 
	\caption{Values of $\delta < 0$ for $3 \leq p \leq 15413$} 
	\label{negdelta}
\begin{center}
    \begin{tabular}{||c|c|c|c|c|c|c||c||}
		\hline
		$\log_{p} \abs{\delta}$ & 0--${1/6}$ & ${1/6}$--${1/3}$ & 
		${1/3}$--${1/2}$ & ${1/2}$--${2/3}$ & ${2/3}$--${5/6}$ &
		${5/6}$--$1$  & total	\\
		\hline
		$\#$ of $p$ &  17 & 78 & 316 & 542 & 51 & 0 & 1004\\
		\hline
	\end{tabular}
\end{center}
\end{table}

It is not clear whether the greater number of negative values of
$\delta$ is significant, or a coincidence of this particular data set.
The mean for Table~\ref{posdelta} is $0.4943$ and the mean for
Table~\ref{negdelta} is $0.5050$.  This may reflect the same apparent
bias towards negative values of $\delta$.

\begin{figure}[h!] 
	\caption{Plot of all values of $\abs{\delta}$ for $3 \leq p \leq 15413$} 
	\label{absplot}
\begin{center}
\includegraphics[height=3in]{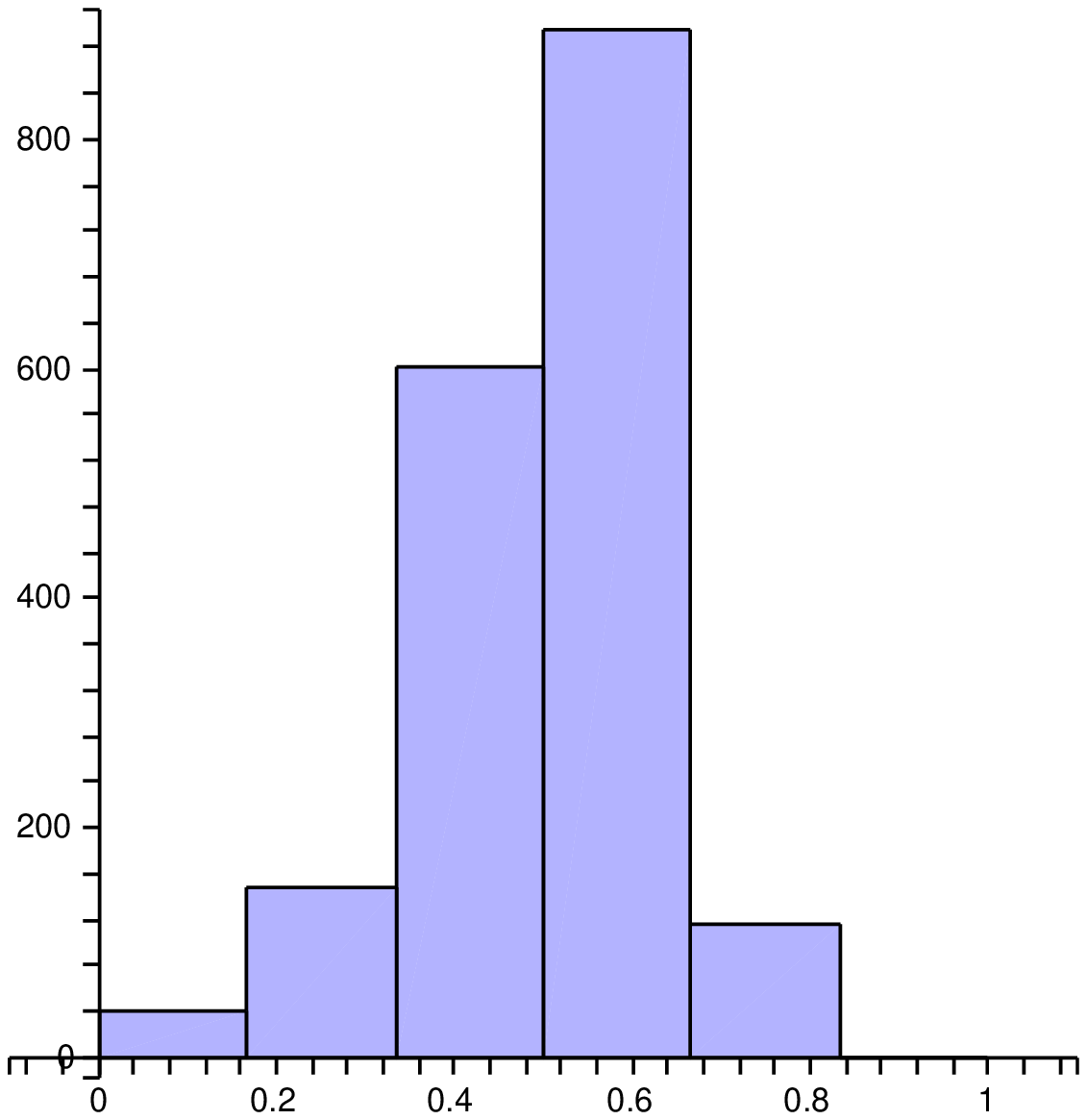}
\end{center}
\end{figure}

Table~\ref{bothtable} and Figure~\ref{absplot} give the values of
$\abs{\delta}$ for all computed values of $\delta$.  The mean for this
table is $0.5003$, which suggests that the expected value of $\log_{p}
\abs{\delta}$ may in fact be $1/2$, i.e., that the values of $\delta$
may cluster around $\sqrt{p}$.  It is not immediately clear how to
derive this from the heuristics.  The sample standard deviation can be
calculated to be $0.1374$, but the data does not appear to be
precisely normally distributed.  This is confirmed by a chi-squared
test for goodness of fit, which returns the extremely small $p$-value
of $7.8039 \cdot 10^{-34}$.\footnote{The $p$-value here can be
interpreted as the chance that a random sample taken from the
predicted distribution would deviate from the distribution as a whole
at least as much as the observed data did.  Thus this set of data is
an extremely bad match for the prediction.  We are using statistical
language in this note even though the data sets do not come from
random variables, and are in fact deterministic.  Thus, all of the
statistical results in this note should be taken with a very large
grain of salt.} A sample skewness of $-0.6785$ and a sample kurtosis
of $3.6516$ can also be computed.  This reflects an asymmetric longer
left tail (toward smaller values of $\log_{p} \abs{\delta}$) and a
somewhat sharper peak than a normal distribution.

\begin{table}[h!] 
	\caption{All values of $\abs{\delta}$ for $3 \leq p \leq 15413$} 
	\label{bothtable}
\begin{center}
    \begin{tabular}{||c|c|c|c|c|c|c||c||}
		\hline
		$\log_{p} \abs{\delta}$ & 0--${1/6}$ & ${1/6}$--${1/3}$ & 
		${1/3}$--${1/2}$ & ${1/2}$--${2/3}$ & ${2/3}$--${5/6}$ &
		${5/6}$--$1$  & total	\\
		\hline
		$\#$ of $p$ &  40 & 147 & 601 & 895 & 116 & 1 & 1800 \\
		\hline
	\end{tabular}
\end{center}
\end{table}

A log-normal distribution was also investigated by taking the
exponential function of the midpoints of each of the class intervals.
This resulted in a mean of $1.6643$, a sample standard deviation of 
$0.2196$, a sample skewness of $-0.2366$ and a sample kurtosis of 
$3.2065$.  Thus the shape of the distribution looks more like a 
normal distribution; however the chi-squared goodness of fit test 
still gives an extremely small $p$-value of $2.2243 \cdot 10^{-10}$.  
Thus this still does not seem to be the correct distribution.  More 
investigation is clearly necessary, both theoretical and statistical.

The data sets from the tables were collected on a Beowulf cluster,
using 16 nodes, each consisting of 2 Pentium~III processors running at
1~Ghz.  The programming was done in C, using MPI, OpenMP, and OpenSSL
libraries.  The collection took approximately 60 hours for the 1800
primes between 3 and 15413 (inclusive).

\section{Conclusion and Future Work}

This note is clearly a preliminary effort.  The fact that we were 
unable to interpret the data as any sort of normal distribution is 
unsatisfying, if not perhaps surprising.  We hope in the future to 
provide at least a conjectural explanation of this data.  A better 
theoretical understanding of the error terms in the theorems we have 
cited would of course be helpful in this.

The project of extending our analysis to three-cycles and more
generally $k$-cycles for small values of $k$, mentioned in previous
papers, still remains to be done.  Along similar lines, Igor
Shparlinski has suggested attempting to analyze the average length of
a cycle.  Daniel Cloutier, a student at the Rose-Hulman Institute of
Technology, has recently begun to collect data which we hope will shed 
light on both of these problems.

\section*{Acknowledgments}

The author would like to thank Pieter Moree for providing the
heuristics for estimating the standard deviation $\sigma$ of $\Nfp_{g
\ANY, h \ANY}(p)$ and for several other results cited in this note.  
He would also like to thank Diane Evans of the Rose-Hulman Institute 
of Technology Mathematics Department for statistical advice.
Finally, he would like to thank the editor for several helpful 
suggestions.

\bibliographystyle{plain}
\bibliography{logcycles}

\end{document}